\documentclass{article}
\usepackage{latexsym}
\usepackage{amssymb}

\begin{document}

\title{An Algebraic Approach \\To Rectangle Packing Problems}

\author{
Baris Altunkaynak\footnote{Department of Physics,
Bogazici University, Bebek, 34342 Istanbul, Turkey.
e-mail: altunkai@boun.edu.tr}
}

\maketitle

\begin{abstract}
A method for converting the geometrical problem of rectangle
packing to an algebraic problem of solving a system of polynomial
equations is described.
\end{abstract}

\newtheorem{Theorem1}{Theorem}
\newtheorem{Proof1}{Proof}

\section{Introduction}
There are many interesting infinite rectangle packing problems
which have been studied. Some of them are: Packing of rectangles
with side lengths $(1/n, 1/(n+1))$ where $n = 1, 2, ...$ into the
unit square, or rectangles of sides $1/n$ where $n = 1, 2, ...$
into a rectangle of area $\pi^2/6$ and so on \cite{Guy, Moser}.
There are many improvements \cite{Chalcraft,Wastlund} for these
problems but some of them are still open. In this paper, we
describe a method to transform a general type of rectangle packing
problem into a
system of polynomial equations. \\

\section{Equivalence of Two Problems}
Let us say we have a set of rectangles (finite or infinite) that we
want to pack into a box (a bigger rectangle) of sides $(A,B)$ which we put
onto a coordinate frame where its left-bottom corner is on the
origin. Let us denote the coordinates of the corners of $n^{th}$
rectangle with $\{(x_{n}^{-}, y_{n}^{-}),(x_{n}^{+},
y_{n}^{-}),(x_{n}^{-}, y_{n}^{+}),(x_{n}^{+}, y_{n}^{+})\}$ as in
Figure 1.

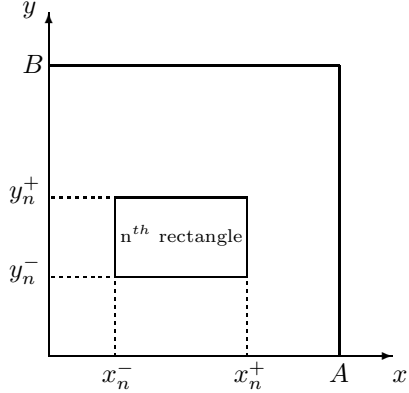
\begin{figure}[htbp]
\begin{center}
\begin{picture}(160,160)
\put(30,10){\vector(0,1){130}} \put(30,10){\vector(1,0){130}}      
\put(30,120){\line(1,0){110}} \put(140,120){\line(0,-1){110}}      

\put(55,40){\line(1,0){50}} \put(105,40){\line(0,1){30}}           
\put(105,70){\line(-1,0){50}} \put(55,70){\line(0,-1){30}}

\put(57,53){\scriptsize{n$^{th}$ rectangle}}                        

\multiput(55,10)(0,3){15} {\line(0,1){1}}                          
\multiput(105,10)(0,3){15} {\line(0,1){1}}                         
\multiput(80,40)(-3,0){17} {\line(-1,0){1}}                        
\multiput(80,70)(-3,0){17} {\line(-1,0){1}}                        

\put(50,0){$x_{n}^{-}$}                           
\put(100,0){$x_{n}^{+}$}                           
\put(15,40){$y_{n}^{-}$}                           
\put(15,70){$y_{n}^{+}$}                           

\put(136,0){$A$} \put(20,117){$B$}                                 
\put(160,0){$x$} \put(20,140){$y$}



\end{picture}
\end{center}
\caption{Coordinates of the corners of the $n^{th}$ rectangle.}
\end{figure}

If we have a perfect packing then we must have the following
equality:

\begin{equation}
\sum_{n} \int_{ y_{n}^{-} }^{ y_{n}^{+} } \!\!\int_{ x_{n}^{-} }^{
x_{n}^{+} } f(x,y)dxdy = \int_{0}^{B} \!\!\!\int_{0}^{A}
f(x,y)dxdy
\end{equation}
for any function $f(x,y)$, because a perfect packing is simply a
covering of the box. This is a necessity condition for
these rectangles to perfectly fit, but it is not
obvious whether it is sufficient or not.

This necessity equation gives us some interesting equalities for
some well known infinite packing problems. For example the problem of
packing the rectangles with sides $(1/n,
1/(n+1))$ with  $n = 1, 2, ...$ into the unit square. For this
problem, if we let $f(x,y)$ to be linear and quadratic functions
of x and y and define $x_{n} = (x_{n}^{+}+x_{n}^{-})/2$,
$y_{n} = (y_{n}^{+}+y_{n}^{-})/2$ (center of mass
coordinates) then we get these following relations:

\begin{eqnarray}
\sum_{n \geq 1} \frac{x_{n}}{n(n+1)} &=& \frac{1}{2}, \\
\sum_{n \geq 1} \frac{y_{n}}{n(n+1)} &=& \frac{1}{2}, \\
\sum_{n \geq 1} \frac{x_{n} y_{n}}{n(n+1)} &=& \frac{1}{4}, \\
\sum_{n \geq 1} \frac{x_{n}^{2} + y_{n}^{2}}{n(n+1)} &=& \frac{1}{3}
+ \frac{\pi^{2}}{36},\\
\sum_{n \geq 1} \frac{(x_{n} + y_{n})^{2}}{n(n+1)} &=& \frac{5}{6}
+ \frac{\pi^{2}}{36},\\
\sum_{n \geq 1} \frac{(x_{n} - y_{n})^{2}}{n(n+1)} &=&
\frac{\pi^{2}}{36} - \frac{1}{6}.
\end{eqnarray}

So a perfect packing must satisfy these equalities. Now by setting
$f(x,y)$ more cleverly we can get a sufficient condition for these
rectangles to perfectly pack.\\

\begin{Theorem1}
Rectangles of sides $(w(n),l(n))$, where $n=1,2,..$, can perfectly
pack a rectangle of sides $(A,B)$ if and only if the following
system of polynomial equations:

\[
\sum_{n} \Big\{(x_{n}^{+})^{S_1}-(x_{n}^{-})^{S_1}\Big\}
\Big\{(y_{n}^{+})^{S_2}-(y_{n}^{-})^{S_2}\Big\}=A^{S_1}B^{S_2} \qquad
S_1,S_2=1,2,...
\]
has a solution with the following constraints:
\begin{eqnarray}
\Delta x_{n}+\Delta y_{n} &=& w(n)+l(n), \nonumber \\
\Delta x_{n} \Delta y_{n} &=& w(n) l(n) \nonumber
\end{eqnarray}
where $\Delta x_{n} = x_{n}^{+}-x_{n}^{-}$ and $\Delta y_{n} = y_{n}^{+}-y_{n}^{-}$. \\
\end{Theorem1}

\begin{Proof1}
In the first equality, let $f(x,y)=e^{px + qy}$, so we get:

\begin{equation}
\sum_{n}
\int_{ y_{n}^{-} }^{ y_{n}^{+} }\!\!\int_{ x_{n}^{-} }^{x_{n}^{+} } e^{px + qy}dxdy
= \int_{0}^{B}\!\!\!\int_{0}^{A} e^{px +qy}dxdy
\qquad \forall p,q \in \mathbb{R}.
\end{equation}

If there exist a perfect packing then we must have the equality above. Now multiply both
sides with $e^{-(pa+qb)}$ where $a, b \in \mathbb{R}$, replace $p$
\& $q$ with $ip$ \& $iq$ and integrate with respect to p and q
from $-\infty$ to $\infty$ to get:
\begin{equation}
\sum_{n} \int_{ y_{n}^{-} }^{ y_{n}^{+} }\!\! \int_{ x_{n}^{-} }^{
x_{n}^{+} } \delta (x-a) \delta (y-b) dxdy = \int_{0}^{B}
\!\!\!\int_{0}^{A} \delta (x-a) \delta (y-b) dxdy.
\end{equation}

One can easily see that the right hand side is equal to 1 if the
point $(a, b)$ is inside the box and 0 if it is outside.
But each of the integrals inside the summation on the left hand
side can only give 0 or 1, so if the right hand side is 0 then
none of the rectangles can contain the point $(a, b)$ which means
that all the rectangles are inside the box and if it is
1 then only one rectangle may contain the point $(a, b)$ which
means that the rectangles do not overlap. This is true for all
points $(a, b)$, hence we must have a perfect packing.

Now if we take the integrals in (8) and cancel some common terms we get:
\begin{eqnarray}
\sum_{n} \left\{ e^{p x_{n}^{+} + q y_{n}^{+}} - e^{p x_{n}^{+} + q
y_{n}^{-}}-e^{p x_{n}^{-} + q y_{n}^{+}} + e^{p x_{n}^{-} + q
y_{n}^{-}} \right\} = \nonumber \\
e^{p A + q B}- e^{p A} - e^{q B} + 1.
\end{eqnarray}
Let's expand this into a series with respect to $p$ and $q$ to
get:
{\setlength\arraycolsep{0pt}
\begin{eqnarray}
\sum_{n} \sum_{k \geq 0} \frac{1}{k!}\Big\{(p x_{n}^{+} + q
y_{n}^{+})^k
&&- (p x_{n}^{+} + q y_{n}^{-})^k - (p x_{n}^{-} + q
y_{n}^{+})^k + (p x_{n}^{-} + q y_{n}^{-})^k \Big\}\nonumber\\
&&= \sum_{k \geq 0} \frac{1}{k!} \Big\{ (p A + q B)^k - (p A)^k - (q
B)^k \Big\} + 1,
\end{eqnarray}
}
and this is equal to:
{\setlength\arraycolsep{0pt}
\begin{eqnarray}
\sum_{n} \sum_{k \geq 0} \sum_{r = 0}^{k}
\frac{1}{k!}
\Big(
\begin{array}{c} k \\ r
\end{array}
\Big)
p^r q^{k-r}
\Big\{&&(x_n^+)^r (y_n^+)^{k-r}-(x_n^+)^r (y_n^-)^{k-r} \nonumber
\\&&- (x_n^-)^r (y_n^+)^{k-r} + (x_n^-)^r (y_n^-)^{k-r}  \Big\} \nonumber \\
= \sum_{k \geq 0} \frac{1}{k!} \Big\{ \sum_{r = 0}^{k}
\Big(
\begin{array}{c} k \\ r
\end{array}
\Big)&&
p^r q^{k-r}A^r B^{k-r}- (p A)^k - (q B)^k \Big\}+1.
\end{eqnarray}
}
Now, $k=0$ term cancel the $+1$ 0n the right hand side, and $r=0$ and $r=k$ terms cancel out
from both sides. If we collect similar terms we get the following equation:
{\setlength\arraycolsep{1pt}
\begin{eqnarray}
\sum_{n} \sum_{k \geq 2} \sum_{r = 1}^{k-1}
\frac{1}{k!}
\Big(
\begin{array}{c} k \\ r
\end{array}
\Big)
&& p^r q^{k-r}
\Big\{ (x_n^+)^r - (x_n^-)^r \Big\} \Big\{(y_n^+)^{k-r}- (y_n^-)^{k-r}\Big\} \nonumber \\
&& = \sum_{k \geq 2}  \sum_{r = 1}^{k-1}\frac{1}{k!}
\Big(
\begin{array}{c} k \\ r
\end{array}
\Big)
p^r q^{k-r}A^r B^{k-r}.
\end{eqnarray}
}
Here $p$ and $q$ are continuous variables, so the above equality holds if and only if
for both sides of the equality the coefficients of the product $p^{S_1} q^{S_2}$
are equal for all integers $S_1, S_2 \geq 1$.
Hence at the end we see that we have these equations:
\begin{equation}
\sum_{n} \Big\{(x_{n}^{+})^{S_1}-(x_{n}^{-})^{S_1}\Big\}
\Big\{(y_{n}^{+})^{S_2}-(y_{n}^{-})^{S_2}\Big\}=A^{S_1}B^{S_2} \qquad
S_1,S_2=1,2,... \nonumber
\end{equation}
Without the constraints the above equation give all possible covering
of the box, one must use the constraints to obtain the desired covering.
Hence the above equations with the constraints are equivalent to the packing problem.
\flushright $\square$ \\
\end{Proof1}

This equation has also a very simple geometrical meaning. Observe
that when we take the product inside the summation we get 4 terms
corresponding to 4 corners of the rectangles. The left-bottom and
right-top corners come with a plus sign and left-top and
right-bottom corners come with a minus sign. For a perfect
packing, rectangles are in contact and for each corner except the
corners of the box there exist another corner which is on the same
point. The terms corresponding to these corners come with opposite
signs and cancel each other. Thus the only remaining terms may
come from the corners of the box, and 3 of them are 0 because of
the corresponding x and y coordinates, and it remains only the
term coming from the right-top corner and this one is exactly
equal to right hand side of the above equation. \\

\section{Discussion}
By using numerical methods, one can solve for a given set of
rectangles their positions corresponding to a perfect packing (if
it exists of course, otherwise we get an empty set) by using a
finite set of equations. This is not an easy task in general but
it can be done as one can try. This numerical solution is not
efficient and impossible for an infinite set of rectangles but the
point is that one can show that a polynomial equation has a root
without explicitly finding it. So we may show that there exist a
perfect packing by showing that these polynomials have a zero (or
vice-versa). A multi-dimensional version of the intermediate value
theorem can be used for this purpose.

\end{document}